%% file: JCAM.tex
\documentclass{elsart}

\usepackage{amsmath,amssymb}
\usepackage{theorem}
\theorembodyfont{\rm}
\usepackage{graphicx}
\usepackage{psfrag}

\newtheorem{problem}{Problem}
\newtheorem{theorem}{Theorem}
 
 \newtheorem{thcor}{Corollary}

\newcommand{\Gd}{{\Gamma_{\rm d}}}
\newcommand{\Gid}{{\Gamma_{\rm id}}}
\newcommand{\Rd}{{R_{\rm d}}}
\newcommand{\Rid}{{R_{\rm id}}}
\newcommand{\bu}{{\overline u}}
\newcommand{\bq}{{\overline q}}
\newcommand{\omegai}{{\omega^*}}
\newcommand{\domega}{{\delta\omega}}
\newcommand{\hu}{{\widehat v}}
\newcommand{\he}{{\widehat e}}

\renewcommand{\epsilon}{\varepsilon}
\renewcommand{\hat}{\widehat}
\newcommand{\uex}{{v^*}}
\newcommand{\huex}{{\hu^*}}
\newenvironment{proof}{{\bf Proof}\quad}%
{\hspace*{\fill}$\square$\medskip}

\begin{document}

\begin{frontmatter}
\title{%
Mathematical Aspects and Numerical Computations of 
an Inverse Boundary Value Identification
Using the Adjoint Method}

\author{Takemi Shigeta}

\address{
Department of Civil Engineering,
National Taiwan University\\
No. 1, Sec. 4, Roosevelt Road, Taipei 10617, Taiwan\\
E-mail: shigeta@ntu.edu.tw}

\begin{abstract}
The purpose of this study is to show some mathematical aspects
of the adjoint method that is a numerical method
for the Cauchy problem, an inverse boundary value problem.
The adjoint method is an iterative method 
based on the variational formulation,
and the steepest descent method minimizes an objective functional
derived from our original problem.
The conventional adjoint method is time-consuming
 in numerical computations because of
the Armijo criterion, which is used to numerically
determine the step size of the steepest descent method.
It is important
to find explicit conditions for the convergence
and the optimal step size.
Some theoretical results about the convergence for the numerical method
are obtained.
Through numerical experiments,
it is concluded that
our theories are effective.
\end{abstract}

\end{frontmatter}

\section{Introduction}
The Cauchy problem is known as an inverse problem.
This problem is to identify unknown boundary value
on a part of the boundary of a bounded domain
for the given boundary data on the rest of the boundary.
In this sense, the problem is regarded as
an inverse boundary value problem.

A numerical method for solution of the Cauchy problem,
proposed by Onishi {\it et al.}\ \cite{Onishi},
is based on the variational formulation.
Namely, this method is constructed by formulating the original problem
to a minimization problem of a functional.
The steepest descent method numerically minimizes the functional.
Finally, the numerical method reduces to an iterative process
in which two direct problems are alternately solved.
The numerical method is also called the adjoint method,
and has been applied to some inverse problems.
The effectiveness of the method has been shown 
by numerical examples \cite{Iijima}, \cite{Shirota}.

However, mathematical properties of the method have not been
clear, yet.
Actually, although the step size of the steepest descent method
 is a very important factor
to influence convergence speed,
it is difficult to obtain theoretically
 the suitable step size in general.
Hence, in the conventional adjoint method,
the step size is numerically determined by the Armijo criterion
\cite{Tosaka}.
But, the Armijo criterion requires many evaluations
because we have to solve the direct problem many times.
It is a big disadvantage that the conventional method is
time-consuming.
Moreover, the convergence of this method has not been proved, yet.

In this paper, we will consider an annulus domain
for simplicity
to prove
that the estimated boundary value obtained by the adjoint method
converges to the exact one.
Moreover, we will
 obtain the suitable step size 
such that the convergence becomes more rapid.
We will confirm the effectiveness of some theoretical results
through simple numerical experiments using the finite element method.

\section{Problem Setting}
For a two dimensional annulus domain
$\Omega:=\{(x,y);\ \Rid^2<x^2+y^2<\Rd^2\}$
with the outer boundary
$\Gd:=\{(x,y);\ x^2+y^2=\Rd^2\}$
and the inner one
$\Gid:=\{(x,y);\ x^2+y^2=\Rid^2\}$,
we consider the Cauchy problem of the Laplace equation:
\begin{problem}\label{prob:cauchy}
For given Cauchy data $(\bu,\bq)\in H^{1/2}(\Gd)\times 
\{\partial v/\partial n\in H^{-1/2}(\Gd);$
$\text{$-\Delta v=0$ in $\Omega$,\
$v|_\Gd=\bu$}$,\ $v|_\Gid=\omega$,\ $\omega\in H^{1/2}(\Gid)\}$,
find $u\in H^{1/2}(\Gid)$ 
 such that
\[
\left\{
\begin{split}
-\Delta u&=0 
&& 
\text{in} 
&&
\Omega,\\
u=\overline u,\quad \frac{\partial u}{\partial n}&=\overline q
&& \text{on} && \Gd,
\end{split}
\right. \]
\end{problem}
where $n$ denotes the unit outward normal to $\Gd$.

\begin{figure}[h]
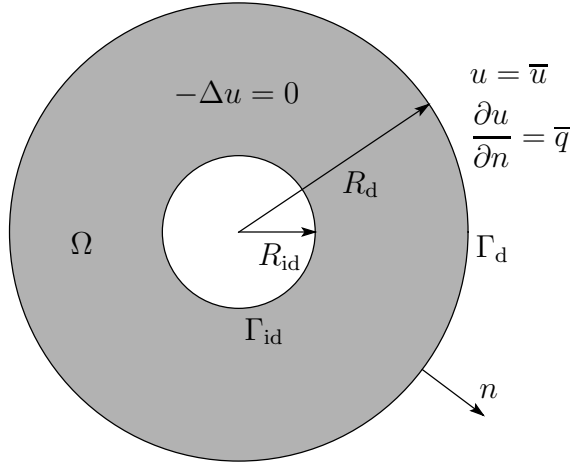

\begin{center}
\input domain.tex
\end{center}
\caption{Cauchy problem}\label{fig:dom}
\end{figure}

In engineering and medical science,
this is an important problem
known as a mathematical model of the inverse problem
of the electrocardiography, which is a problem to identify
unknown electric potential on the surface of a human heart
by observing the electric potential on the surface of a human body.

\section{Variational Formulation}
To solve Problem \ref{prob:cauchy}, we consider the following
variational problem based on the method of the least squares:
\begin{problem}
Find $\omegai\in H^{1/2}(\Gid)$ such that
\begin{equation}
 J(\omegai)=\inf_{\omega\in H^{1/2}(\Gid)} J(\omega), 
 \label{eq:functional}
\end{equation}
where the functional $J$ is defined as
\[ J(\omega):=\int_\Gd|v(\omega)-\bu|^2d\Gamma, \]
and $v=v(\omega)\in H^1(\Omega)$ depending on the given boundary value
$\omega\in H^{1/2}(\Gid)$
is the solution of the following mixed boundary value problem
called the primary problem:
\begin{equation} \left\{
\begin{split}
-\Delta v&=0 
&& 
\text{in} &&
\Omega,\\
\frac{\partial v}{\partial n}&=\bq
&& \text{on} && \Gd,\\
v&=\omega
&& \text{on} && \Gid.
\end{split} \right. \label{eq:primary}
\end{equation}
\end{problem}

Our strategy is to find the minimum of $J$
by generating a minimizing sequence $\{\omega_{k}\}_{k=0}^\infty
\subset H^{1/2}(\Gid)$
by the steepest descent method:
\begin{equation}
 \omega_{k+1}=\omega_{k}-\rho_kJ'(\omega_{k}),
\label{eq:sdm}
\end{equation}
starting with an initial guess $\omega_{0}\in H^{1/2}(\Gid)$,
with a suitably chosen step size $\{\rho_k\}_{k=0}^\infty$.
The derivative $J'(\omega)$ is the first variation of $J$,
defined by
\[ J(\omega+\domega)-J(\omega)=( J'(\omega),\domega )
+o(\|\domega\|), \]
where $(\cdot,\cdot)$ and $\|\cdot\|$
denote the inner product
and the norm in $L^2(\Gid)$, respectively:
\[ ( f,g ):=\int_\Gid fg\,d\Gamma,
\qquad
\|f\|:=( f,f)^{1/2}. \]

The first variation is explicitly given \cite{Iijima}, \cite{Onishi} --
\cite{Tosaka} by
\begin{equation}
\left.
J'(\omega)=
-\frac{\partial\hu}{\partial n}\right|_\Gid,
\label{eq:fv}
\end{equation}
where $\hu=\hu(v(\omega))\in H^2(\Omega)$ depending on
the solution $v=v(\omega)$ of the primary problem 
is the solution of the following mixed boundary value problem
called the adjoint problem:
\begin{equation}
 \left\{
\begin{split}
-\Delta \hu&=0 
&&
\text{in} 
&&
\Omega,\\
\frac{\partial \hu}{\partial n}&=2(v-\bu)
&& \text{on} && \Gd,\\
\hu&=0
&& \text{on} && \Gid.
\end{split}
\right. \label{eq:adjoint}
\end{equation}
We remark that $J'(\omega)\in H^{1/2}(\Gid)$.

\section{Algorithm}
Our numerical method for the Cauchy problem
reduces to an iterative process to minimize the functional
 (\ref{eq:functional})
by the steepest descent method (\ref{eq:sdm})
after solving the primary problem (\ref{eq:primary})
and the adjoint problem (\ref{eq:adjoint})
to obtain the first variation (\ref{eq:fv}).

For given $\omega=\omega_k$,
we denote the solutions of the primary and the adjoint problems
by $v_k:=v(\omega_k)$, $\hu_k:=\hu(v(\omega_k))$, respectively.
Then, our algorithm can be summarized as follows:

\emph{Algorithm}
\begin{description}
\item[Step 0.] Pick an initial guess $\omega_{0}\in H^{1/2}(\Gid)$,
and set $k:=0$.
\item[Step 1.] Solve the primary problem
\begin{equation}
\left\{
\begin{split}
-\Delta v_k&=0 &  & \text{in} &  & \Omega, \\
\frac{\partial v_k}{\partial n}&=\bq
&& \text{on} && \Gd,  \\
v_k&=\omega_k
&& \text{on} && \Gid
\end{split}
\right. \label{eq:prim}
\end{equation}
to find $v_k|_\Gd\in H^{1/2}(\Gd)$.

\item[Step 2.] Solve the adjoint problem
\begin{equation}
\left\{
\begin{split}
-\Delta \hu_k&=0 &  & \text{in} &  & \Omega,  \\ 
\frac{\partial \hu_k}{\partial n}&=2(v_k-\bu)
&& \text{on} && \Gd, \\ 
\hu_k&=0
&& \text{on} && \Gid  
\end{split}
\right. \label{eq:adj}
\end{equation}
to find the first variation
\[ \left. J'(\omega_k)=
-\frac{\partial\hu_k}{\partial n}
\right|_\Gid
\in H^{1/2}(\Gid).
 \]
\item[Step 3.] Choose the step size $\rho_k$.
\item[Step 4.] 
Update the boundary value by
the steepest descent method
\begin{equation}
 \omega_{k+1}=\omega_k-\rho_kJ'(\omega_k). 
\label{eq:sdm2}
\end{equation}
\item[Step 5.] Set $k:=k+1$, and go to Step 1.
\end{description}

In numerical computations, the primary and the adjoint problems
in Steps 1 and 2 are numerically solved by
the triangular finite element method.
The conventional method to choose the suitable step size $\rho_k$ in Step 3
is the Armijo criterion,
which guarantees for the sequence $\{\rho_k\}_{k=0}^\infty$
to satisfy
\begin{equation}
 J(\omega_k-\rho_kJ'(\omega_k))\leq
J(\omega_k)-\xi\rho_k\|J'(\omega_k)\|^2 
\label{eq:armijo}
\end{equation}
with a constant $0<\xi<1/2$.

{\samepage
\emph{Controlling the step size}
\begin{description}
\item[Step 3.0.] Give $0<\xi<1/2$, $0<\tau<1$ and
the sufficiently small $\epsilon>0$.
\item[Step 3.1.] If $\|J'(\omega_{k})\|<\epsilon$, then stop;
else go to Step 3.2.
\item[Step 3.2.] Set $\beta_0:=1$, $m:=0$.
\item[Step 3.3.] If $J(\omega_{k}-\beta_mJ'(\omega_{k}))
\leq J(\omega_{k})-\xi\beta_m\|J'(\omega_{k})\|^2$,
then set $\rho_k:=\beta_m$; else go to Step 3.4.
\item[Step 3.4.] Set $\beta_{m+1}:=\tau\beta_m$.
\item[Step 3.5.] Set $m:=m+1$, and go to Step 3.3.
\end{description}
}

To choose the step size, we have to evaluate
the functional value on the left hand side in (\ref{eq:armijo}) many times.
It means that the primary problem has to be solved many times.
This is an disadvantage of the conventional method.

\section{Convergence Proof and the Optimal Step Size}
In this section, we prove that our numerical method is convergent
under some assumption.
Moreover, we propose suitable step sizes
to avoid the disadvantage of the conventional method.
The following argument is based on the convergence proof of
the Dirichlet-Neumann alternating method \cite{Yu2},
which is one of the domain decomposition methods \cite{Lu}.

We denote by $\uex$ and $\huex$
the solutions of the primary and the adjoint problems
for the boundary value $\omega=\omegai$, respectively.
Let the error functions $e_k:=\uex-v_k$, $\he_k:=\huex-\hu_k$,
and $\mu_k:=\omegai-\omega_k$.
Then, from (\ref{eq:prim}) and (\ref{eq:adj}), 
we see that the functions $e_k$ and $\he_k$ are the solutions of
\begin{equation}
\left\{
\begin{split}
-\Delta e_k&=0 &  & \text{in} &  & \Omega,\\
\frac{\partial e_k}{\partial n}&=0
&& \text{on} && \Gd,\\
e_k&=\mu_k
&& \text{on} && \Gid,
\end{split}
\right. \label{eq:prim2}
\end{equation}
and
\begin{equation}
 \left\{
\begin{split}
-\Delta \he_k&=0 &  & \text{in} &  & \Omega,\\
\frac{\partial \he_k}{\partial n}&=2e_k
&& \text{on} && \Gd,\\
\he_k&=0
&& \text{on} && \Gid,
\end{split}
\right. \label{eq:adj2}
\end{equation}
respectively.

We assume that the error $\mu_k$ can be expanded
into the finite Fourier series:
\begin{equation}
 \mu_k=\sum_{|j|=M}^Na_j^{(k)}e^{ij\theta}
\label{eq:mu}
\end{equation}
with some integers $M$ and $N$ such that $N\geq M\geq 0$.
From (\ref{eq:prim2}), we have
\begin{equation}
 e_k|_\Gd=\sum_{|j|=M}^N
\frac{2\Rd^{|j|}\Rid^{|j|}}
{\Rd^{2|j|}+\Rid^{2|j|}}a_j^{(k)}
e^{ij\theta}.
\label{eq:ek}
\end{equation}
Substituting (\ref{eq:ek}) into (\ref{eq:adj2}),
we can derive
\begin{equation}
 \left.\frac{\partial \hat e_k}{\partial n}\right|_\Gid
=-\sum_{|j|=M}^N
\frac{8\Rd^{2|j|+1}\Rid^{2|j|-1}}{(\Rid^{2|j|}+\Rd^{2|j|})^2}
a_j
e^{ij\theta}.
\label{eq:en}
\end{equation} 
Here, noting that $\huex=0$, we have
\begin{equation}
\left. J'(\omega_k) 
=\frac{\partial\hat e_k}{\partial n}\right|_\Gid .
\label{eq:j1}
\end{equation}
We can see from (\ref{eq:sdm2}) that
\begin{equation}
 \mu_{k+1}=\mu_k+\rho_kJ'(\omega_k).
\label{eq:mu_sdm}
\end{equation}
Substituting (\ref{eq:mu}) and (\ref{eq:j1}) (namely (\ref{eq:en}))
into (\ref{eq:mu_sdm}), we have
\begin{equation} a_j^{(k+1)}
=(1-C_j\rho_k)a_j^{(k)} \label{eq:ajk1}
\end{equation}
with
\[ C_j:=\frac{8\Rd^{2|j|+1}\Rid^{2|j|-1}}{(\Rid^{2|j|}+\Rd^{2|j|})^2}. \]
If we take the step size as
\[ 0<\rho_k<
\frac 2{C_M},
\]
then we have
\begin{equation}
 |\delta_j^{(k)}|\leq\delta^{(k)}
=\max\left\{
\left|
1-C_M\rho_k
\right|,
\left|
1-C_N\rho_k
\right|
\right\}<1
\label{eq:delta}
\end{equation}
with
$\delta_j^{(k)}:=1-C_j\rho_k$ and
$\displaystyle\delta^{(k)}:=\max_{M\leq j\leq N}|\delta_j^{(k)}|$.
It follows that $\{\|\mu_k\|\}_{k=0}^\infty$ is a strictly monotone decreasing
sequence:
\[ \|\mu_{k+1}\|\leq \delta^{(k)} \|\mu_k\| \]
with the compression factor $\delta^{(k)}<1$.
Therefore, we obtain that $\mu_k\to 0$, that is,
$\omega_k\to\omegai$ as $k\to\infty$.

We notice that the convergence cannot be guaranteed in general
because $\delta^{(k)}=1$ in (\ref{eq:delta})
if $\mu_k$ is expanded into the infinite series.

As a consequence, we can obtain the following theorem:
\begin{theorem}
For the exact boundary value $\omegai$,
we assume that the error $\mu_k=\omegai-\omega_k$ 
can be expanded into the finite Fourier series:
\begin{equation}
 \mu_k=\sum_{|j|=M}^Na_j^{(k)}e^{ij\theta}
\label{eq:mu2}
\end{equation}
with some integers $M$ and $N$ such that $N\geq M\geq 0$.
If we take the step size as
\[ 0<\rho_k<
\frac 2{C_M},
\qquad C_j:=\frac{8\Rd^{2|j|+1}\Rid^{2|j|-1}}{(\Rid^{2|j|}+\Rd^{2|j|})^2} ,
\]
then $\{\omega_k\}_{k=0}^\infty$ converges to $\omegai$.
\end{theorem}

\begin{thcor}\label{cor:opt}
The optimal step size $\rho_{\rm opt}$ 
in the sense that the compression factor is minimized is given by
\begin{equation}
\rho_{\rm opt}=\frac 2{C_M+C_N} .
\label{eq:rhoopt}
\end{equation}
Then, the optimal compression factor $\delta_{\rm opt}$ is given by
\begin{equation}
 \delta_{\rm opt}=\frac{C_M-C_N}{C_M+C_N} .
\label{eq:deltaopt}
\end{equation}
\end{thcor}
\begin{proof}
Taking
\[ |1-C_M\rho_k|=|1-C_N\rho_k|=\delta_{\rm opt} \]
so as to minimize  $\delta^{(k)}$ in (\ref{eq:delta}),
we obtain (\ref{eq:rhoopt}) and (\ref{eq:deltaopt}).
\end{proof}

The next theorem is very effective in actual computations.

\begin{theorem}\label{th:2}
The error $\mu_k$ is assumed to be expanded into (\ref{eq:mu2}).
If we put $\rho_k=1/C_{M+k}$ or $\rho_k=1/C_{N-k}$ 
($k=0,1,\ldots,N-M$),
then we have $\mu_{N-M+1}=0$, namely 
the exact $\omega^*$ is obtained after $(N-M+1)$ iterations.
\end{theorem}
\begin{proof}
For $k=0,1,\ldots,N-M$,
substituting $\rho_k=1/C_{M+k}$ (resp. $\rho_k=1/C_{N-k}$)
 into (\ref{eq:ajk1}),
we have $a_{M+k}^{(k+1)}=0$ (resp. $a_{N-k}^{(k+1)}=0$).
Then, it follows that
$a_{M+k}^{(l)}=0$ (resp. $a_{N-k}^{(l)}=0$)
for all $l=k+1,k+2,\ldots,N-M+1$.
\end{proof}

\section{Numerical Experiments}
Let the radii $\Rd=3$ and $\Rid=1$.
The domain $\Omega$ is divided into triangular finite elements
with 8552 elements and 4436 nodes.
We take $\omega_0=0$ as an initial guess.
The stop condition of our calculations is given
by $J(\omega_k)<10^{-5}$.

\begin{figure}[h]
\begin{center}
\includegraphics[width=7cm]{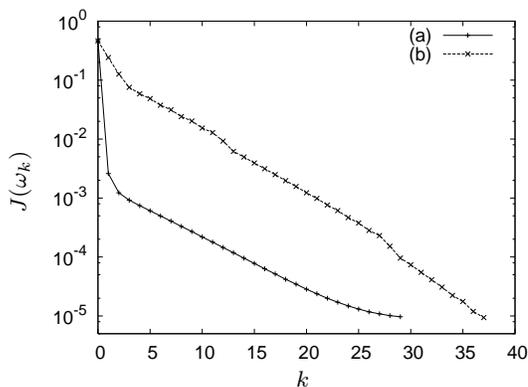}
\caption{Functional values
}\label{fig:sol}
\end{center}
\end{figure}

\begin{figure}
\begin{center}\hspace*{-7em}
\includegraphics[width=10cm]{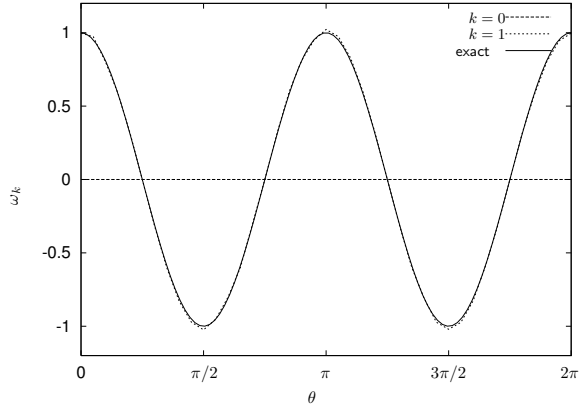}
\end{center}
\caption{The estimated boundary value $\omega_k$}\label{fig:omg1}
\end{figure}

As the first example, let the Cauchy data be
$(\bu,\bq)=(9\cos 2\theta,6\cos 2\theta)$.
Then, the exact potential in $\Omega$
and the related exact boundary value on $\Gid$
can be written as follows:
\[ \left\{
\begin{split}
u^*&=r^2\cos 2\theta, \\
\omega^*&=u^*|_\Gid=\cos 2\theta.
\end{split}
\right. \]
We can see from the Cauchy data that $M=N=2$.
It follows that $\rho_{\rm opt}=1/C_2=1681/486(\approx 3.46)$ and
$\delta_{\rm opt}=0$.
Hence, theoretically it should hold that $\omega_1=\omega^*$.
But, in actual computations, we remark that
$\omega_1\neq\omega^*$ due to discretization errors by
the finite element method.
For $k\geq 1$, we regard that $M=0$ and $N$ is sufficiently large.
Then, according to Corollary \ref{cor:opt}, we have 
$\rho_{\rm opt}\approx 2/C_0=\Rid/\Rd=1/3$.
Hence, we take
\begin{equation}
 \rho_k=\left\{
\begin{array}{ll}
1681/486(\approx 3.46) & (k=0) \\
1/3(\approx 0.33) & (k\geq 1)
\end{array}
\right. .\label{eq:optex}
\end{equation}
Then, the graph (a) in Figure~\ref{fig:sol}
shows the variations of functional value.
Since the functional value decreases steeply
by taking the step sizes as large as possible,
the convergence becomes rapid.
On the other hand, the graph (b)
is the result in the case when we apply
 the Armijo criterion with $\xi=1/3$ and $\tau=1/2$.
When the Armijo criterion is applied,
only to choose the step sizes,
the primary problem has to be solved 89 times.
Hence, until the estimated boundary value converges,
in all we have to solve the direct problems
$58(=29\times 2)$ times in the case of the graph (a)
 and $163(=37\times 2+89)$ times
in the case of the graph (b), respectively.
We can see that (\ref{eq:optex}) makes convergence more rapid
than the Armijo criterion.
On the other hand, Figure~\ref{fig:omg1} shows the estimated
boundary value $\omega_k$ and the exact one $\omega^*$.
We can see that
$\omega_1$ is in good agreement with $\omega^*$.
It is concluded that the estimated boundary value
quickly converges to the exact one.

\begin{figure}[h]
\begin{center}
\includegraphics[width=7cm]{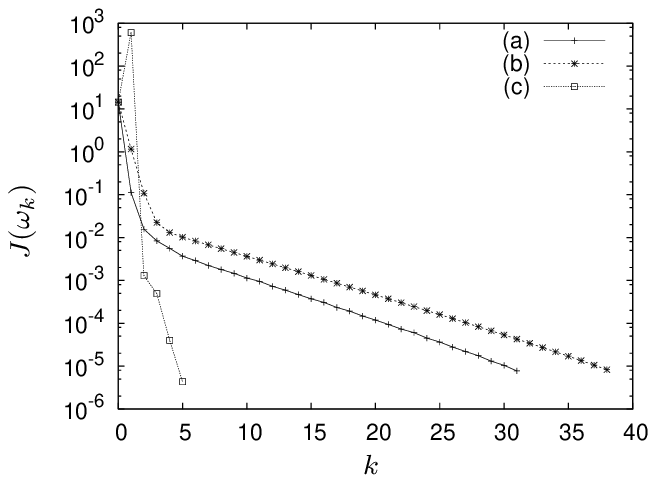}
\caption{Functional values
}\label{fig:sol2}
\end{center}
\end{figure}

As the second example, we assume that the Cauchy data is
\[ (\bu,\bq)=
\left(6\sin \theta-\frac 32\cos\theta+\frac 94\cos 2\theta,
2\sin \theta-\frac 12\cos\theta+\frac 32\cos 2\theta\right). \]
Then, the exact potential in $\Omega$
and the related exact boundary value on $\Gid$
can be written as follows:
\[ \left\{
\begin{split}
u&=r(2\sin\theta-\frac 12\cos\theta)+\frac 14r^2\cos 2\theta, \\
\omega^*&=u|_{\Gamma_{\rm id}}=
2\sin\theta-\frac 12\cos\theta+\frac 14\cos 2\theta .
\end{split}
\right.
\]
We can see from the Cauchy data that $M=1$ and $N=2$.
But, due to numerical errors,
we assume that $M=0$.
According to Theorem \ref{th:2}, we take the step size as
\begin{equation}
 \rho_k=\left\{
\begin{array}{ll}
1/C_{2-k} & (k=0,1,2) \\
1/3(\approx 0.33) & (k\geq 3)
\end{array}
\right. .\label{eq:optex2}
\end{equation}
Then, the graph (c) in Figure~\ref{fig:sol2}
shows the variations of functional value.
As we can see, the functional value reaches less than $10^{-5}$ at $k=5$.
On the other hand, 
using the Armijo criterion with $\xi=1/3$ and $\tau=1/2$,
we obtain the functional value less than $10^{-5}$ at $k=31$,
which is shown as the graph (a).
When we apply the Armijo criterion,
only to choose the step sizes,
the primary problem has to be solved 76 times.
Hence, until the estimated boundary value converges,
in all we have to solve the direct problems
$10(=5\times 2)$ times in the case of the graph (c)
 and $138(=31\times 2+76)$ times in the case of the graph (a), respectively.
If we do not know the concrete values of $M$ and $N$,
it is reasonable to assume that $M=0$ and $N$ is sufficiently large.
Then, the step size is taken as
$\rho_k=1/3(\approx 0.33)$
for all $k\geq 0$ according to Corollary \ref{cor:opt}.
The variations of functional value for this step size
are shown by the graph (b).
Although the number of iterations for the graph (b)
is greater than that for the graph (a),
the computational cost for (b) is less than that for (a) in all.

\begin{figure}
\begin{center}\hspace*{-7em}
\includegraphics[width=10cm]{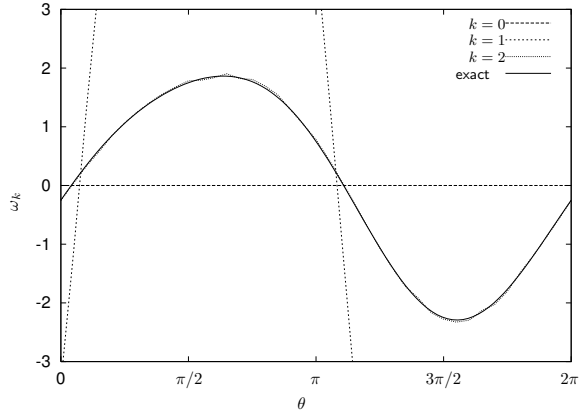}
\end{center}
\caption{The estimated boundary value $\omega_k$}\label{fig:omg2}
\end{figure}

Figure~\ref{fig:omg2} shows the estimated
boundary value $\omega_k$ and the exact one $\omega^*$.
We can see that $\omega_2$ is roughly the same as $\omega^*$.

Therefore, through two examples, it is concluded that
the computational cost can reduce
and the estimated boundary value quickly converges to the exact one
by applying Corollary \ref{cor:opt} and Theorem \ref{th:2},
which are very effective in numerical computations.

\section*{Acknowledgements}
This research was partially supported by
the Ministry of Education, Culture, Sports, Science and Technology,
Grant-in-Aid for
Scientific Research (Young Scientists (B),
No.~18740061),
2006 when the author worked at
Tokyo University of Science, Yamaguchi, Japan. 
The author is most grateful to Professor D.~L.~Young
at National Taiwan University,
who kindly gave helpful suggestions about this paper.

\end{document}

%% file: domain.tex
\unitlength 0.1in
\begin{picture}( 24.8000, 24.0000)(  0.0000,-24.0000)
%
\special{pn 8}%
\special{sh 0.300}%
\special{ar 1200 1200 1200 1200  0.0000000 6.2831853}%
\put(17.4000,-8.8000){\makebox(0,0)[lt]{$\Rd$}}%
\put(24.4000,-12.2000){\makebox(0,0)[lt]{$\Gd$}}%
\put(12.3000,-16.5000){\makebox(0,0)[lt]{$\Gid$}}%
%
\special{pn 8}%
\special{sh 0}%
\special{ar 1200 1200 400 400  0.0000000 6.2831853}%
%
\special{pn 8}%
\special{pa 1198 1202}%
\special{pa 2194 532}%
\special{fp}%
\special{sh 1}%
\special{pa 2194 532}%
\special{pa 2128 554}%
\special{pa 2150 562}%
\special{pa 2150 586}%
\special{pa 2194 532}%
\special{fp}%
\put(13.0000,-12.6000){\makebox(0,0)[lt]{$\Rid$}}%
%
\special{pn 8}%
\special{pa 2160 1920}%
\special{pa 2480 2160}%
\special{fp}%
\special{sh 1}%
\special{pa 2480 2160}%
\special{pa 2440 2104}%
\special{pa 2438 2128}%
\special{pa 2416 2136}%
\special{pa 2480 2160}%
\special{fp}%
\put(8.6000,-4.1000){\makebox(0,0)[lt]{$-\Delta u=0$}}%
\put(3.2000,-12.0000){\makebox(0,0)[lt]{$\Omega$}}%
\put(24.1000,-3.2000){\makebox(0,0)[lt]{$u=\overline{u}$}}%
\put(24.1000,-5.2000){\makebox(0,0)[lt]{$\displaystyle\frac{\partial u}{\partial n}=\overline{q}$}}%
%
\special{pn 8}%
\special{pa 1200 1200}%
\special{pa 1600 1200}%
\special{fp}%
\special{sh 1}%
\special{pa 1600 1200}%
\special{pa 1534 1180}%
\special{pa 1548 1200}%
\special{pa 1534 1220}%
\special{pa 1600 1200}%
\special{fp}%
\put(24.6000,-20.7000){\makebox(0,0)[lb]{$n$}}%
\end{picture}%